\documentstyle[11pt,leqno,here]{article}
\input amssym.def
\input amssym.tex
\input psfig

\newcommand{\beql}[1]{\begin{equation}\label{#1}}
\newcommand{\eeq}{\end{equation}}

\makeatletter
\@addtoreset{figure}{section}
\def\thefigure{\thesection.\@arabic\c@figure}
\@addtoreset{table}{section}
\def\thetable{\thesection.\@arabic\c@table}

\def\@sect#1#2#3#4#5#6[#7]#8{\ifnum #2>\c@secnumdepth
     \def\@svsec{}\else
     \refstepcounter{#1}\edef\@svsec{\csname the#1\endcsname.\hskip .75em }\fi
     \@tempskipa #5\relax
      \ifdim \@tempskipa>\z@
        \begingroup #6\relax
          \@hangfrom{\hskip #3\relax\@svsec}{\interlinepenalty \@M #8\par}%
        \endgroup
       \csname #1mark\endcsname{#7}\addcontentsline
         {toc}{#1}{\ifnum #2>\c@secnumdepth \else
                      \protect\numberline{\csname the#1\endcsname}\fi
                    #7}\else
        \def\@svsechd{#6\hskip #3\@svsec #8\csname #1mark\endcsname
                      {#7}\addcontentsline
                           {toc}{#1}{\ifnum #2>\c@secnumdepth \else
                             \protect\numberline{\csname the#1\endcsname}\fi
                       #7}}\fi
     \@xsect{#5}}
\def\@begintheorem#1#2{\it \trivlist \item[\hskip \labelsep{\bf #1\ #2.}]}
\def\section{\@startsection {section}{1}{\z@}{-3.5ex plus -1ex minus
 -.2ex}{2.3ex plus .2ex}{\normalsize\bf}}

\makeatother

\begin{document}
\begin{center}
{\Large {\bf CURVED HEXAGONAL PACKINGS OF EQUAL DISKS IN A CIRCLE}}\\
\vspace{1\baselineskip}
{\em B. D. Lubachevsky} \\
{\em R. L. Graham} \\
\vspace{.25\baselineskip}
AT\&T Bell Laboratories, \\
Murray Hill, New Jersey 07974 \\
\vspace{1\baselineskip}
\vspace{1.5\baselineskip}
{\bf ABSTRACT}
\vspace{.5\baselineskip}
\end{center}

\setlength{\baselineskip}{1.5\baselineskip}

For each $k \ge 1$ and corresponding 
hexagonal number $h(k) = 3k(k+1)+1$,
we introduce $m(k) = \max \{ \frac {(k-1)!} {2}, 1\}$ packings of
$h(k)$ equal disks inside a circle which 
we call the {\em curved hexagonal} packings.
The curved hexagonal packing of 7 disks ($k = 1$, $m(1)=1$) is well known
and the one of 19 disks ($k = 2$, $m(2)=1$) 
has been previously conjectured to be optimal. 
New curved hexagonal packings of 37, 61, and 91 disks 
($k = 3$, 4, and 5, $m(3)=1$, $m(4)=3$, and $m(5)=12$)
were the densest we obtained on a computer
using a so-called ``billiards'' simulation algorithm.
A curved hexagonal packing pattern is invariant under 
a $60 ^{\circ}$ rotation.
For $k \rightarrow \infty$, the density (covering fraction) 
of curved hexagonal packings tends to $\frac {\pi^2} {12}$.
The limit is smaller than the density of the known optimum disk packing in
the infinite plane.
We found disk configurations that are denser than curved hexagonal packings
for 127, 169, and 217 disks ($k = 6$, 7, and 8).

In addition to new packings for $h(k)$ disks,
we present new packings we found for $h(k)+1$ and $h(k)-1$ disks
for $k$ up to 5, i.e., for 36, 38, 60, 62, 90, and 92 disks.
The additional packings show the ``tightness''
of the curved hexagonal pattern for $k \le 5$:
deleting a disk does not change
the optimum packing and its quality significantly,
but adding a disk causes a substantial rearrangement
in the optimum packing and substantially decreases the quality.
 \\

\vspace{1.\baselineskip}
\clearpage
\large\normalsize
\renewcommand{\baselinestretch}{1}
\thispagestyle{empty}
\setcounter{page}{1}
\begin{center}
{\Large {\bf CURVED HEXAGONAL PACKINGS OF EQUAL DISKS IN A CIRCLE}} \\
\vspace{1\baselineskip}
{\em B. D. Lubachevsky} \\
{\em R. L. Graham} \\
\vspace{.25\baselineskip}
AT\&T Bell Laboratories, \\
Murray Hill, New Jersey 07974 \\
\vspace{1\baselineskip}
\end{center}
\setlength{\baselineskip}{1.46\baselineskip}
\section{Introduction}
\hspace*{\parindent}
Patterns of dense geometrical packings are sensitive
to the geometry of the enclosing region of space.
In particular, dense packings of equal nonoverlapping disks
in a circle are different from
those in a regular hexagon,
as one might expect
(see, e.g., \cite{CFG}, \cite{G} or \cite{R}).
\begin{figure}[htb]
\centerline{\psfig{file=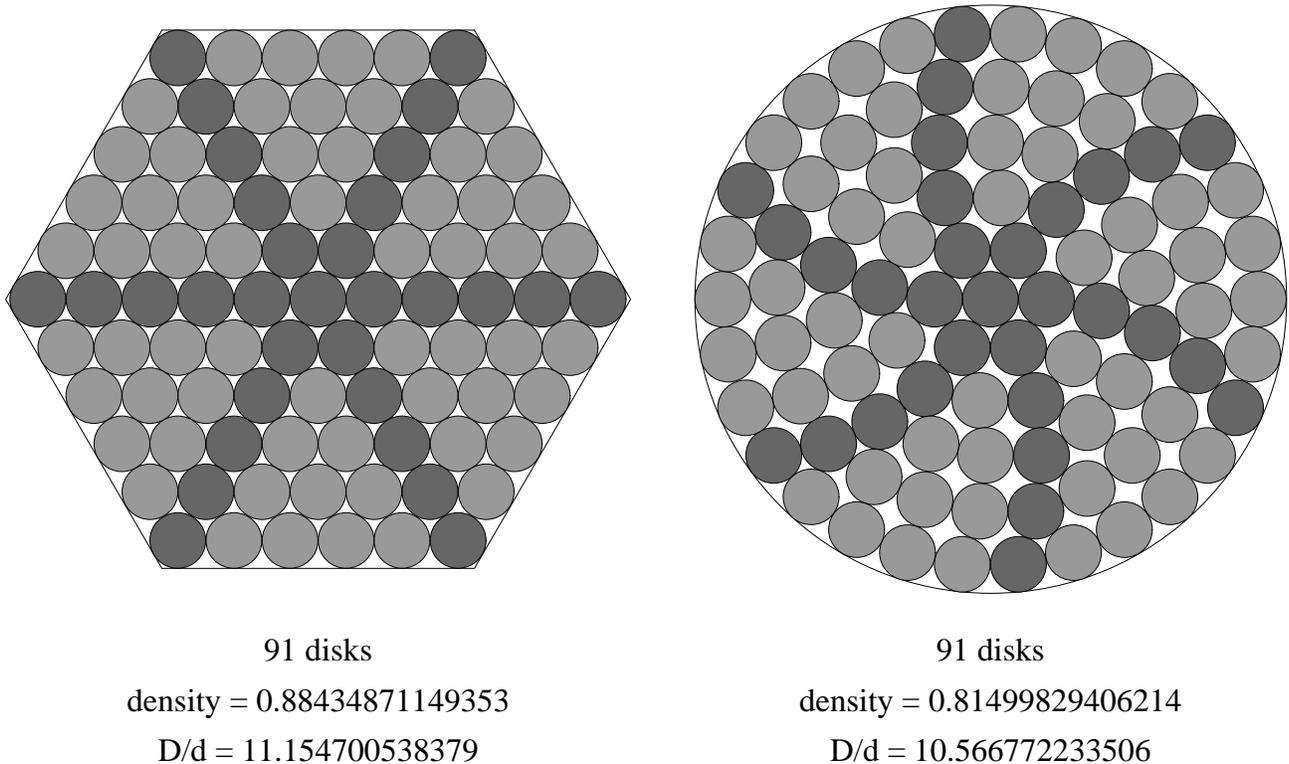,width=7in}}
\caption{{\it Left}. The well known best packing of $h(5)=91$ disks 
in a regular hexagon. 
{\it Right}. One of the 12 best (that we found)
packings of $h(5)=91$ disks in a circle.
In both diagrams,
density shown is the ratio of the area covered by disks
to the area of the container,
$D$ and $d$ are the diameters of the container and of
the disks, respectively.
(The different shading of the disks is auxiliary;
it is not a part of the pattern.)}
\label{twolarge}
\end{figure}
In this paper, for each $k \ge 1$ 
and corresponding hexagonal number $h(k) = 3k(k+1)+1$,
we present a pattern of packings of $h(k)$ equal disks in a circle
which can be viewed a ``curved'' analogue of the 
densest packing of $h(k)$ disks in a regular hexagon.
For a particular $k$, there exists a set of
$m(k) = \max \{ \frac {(k-1)!} {2}, 1 \}$ different
{\em curved hexagonal packings} of the same quality.
A curved hexagonal packing pattern is invariant under
a $60 ^{\circ}$ rotation.
The density (covering fraction) 
of a curved hexagonal packing
tends to $\frac {\pi^2} {12}$ as $k \rightarrow \infty$.
Because the limit is smaller than
the density of the best
(hexagonal) packing of equal disks on an infinite plane,
the curved hexagonal packing of $h(k)$ disks
can not be optimal for a sufficiently large $k$.

It is remarkable, though, that for several initial
values of $k$ there seems to be no better packing 
than the curved hexagonal ones.
Indeed, for 7 disks ($k=1$, $m(1)=1$) the curved hexagonal packing
is well known to be optimal
and the one for 19 disks ($k=2$, $m(2)=1$) has been previously
conjectured as such \cite{K}.
For 37, 61, and 91 disks ($k=3$, 4, and 5, $m(3)=1$, $m(4)=3$, $m(5)=12$),
the curved hexagonal packings were the densest we obtained by 
computer experiments
using the so-called ``billiards'' simulation algorithm.

The ``billiards'' simulation algorithm \cite{L} \cite{LS}
has so far proved to be a reliable method for generating optimal
packings of disks in an equilateral triangle \cite{GL1}.
Our experiments with this algorithm for packings in a circle
either confirmed or improved the best previously reported packings
for $n \le 25$. 
We are unaware of any published conjectures for packing $n > 25$
disks in a circle,
but the ``billiards'' algorithm kept producing packings for many $n > 25$,
specifically, for $n = h(3) = 37$, $n = h(4) = 61$, and $n = h(5) = 91$.
(A detailed account of these experiments merged with the experiments
of Nurmela and \"Osterg\aa rd will be reported in a forthcoming paper
\cite{GLNO}.)
The latter three sets of packings happened to have the curved hexagonal pattern
and they were the best found for their value of $n$.
As for $n = h(6) = 127$, $n = h(7) = 169$, and $n = h(8) = 217$, 
the algorithm found packings which are
{\em better} than the corresponding curved hexagonal packings.

For the values $k \le 5$, for which the densest packings of $h(k)$
disks in a circle apparently have the curved hexagonal patterns,
these packings look ``tight.''
To test our intuition of their ``tightness''
we compared these packings with packings
obtained for $h(k)-1$ and $h(k)+1$ disks.
Thus, we generated dense packings for 36, 38, 60, 62, 90, and 92 disks
and verified that deleting a disk does not change the
optimum packing and its quality significantly,
but adding a disk causes a substantial rearrangement in the optimum packing
and substantially decreases the quality.
This tightness may be considered an analogue of
the similar tightness property for the infinite classes
of packings in an equilateral triangle
as noted in \cite{GL1}. 
Specifically the variations in the packing pattern and quality
when one disk is added or subtracted
are similar to those observed for
packings of $\frac {k(k+1)} {2}$ disks
in the triangles.
\section{Packings of 7, 19, 37, and 61 disks in a circle}
\begin{figure}[ht]
\centerline{\psfig{file=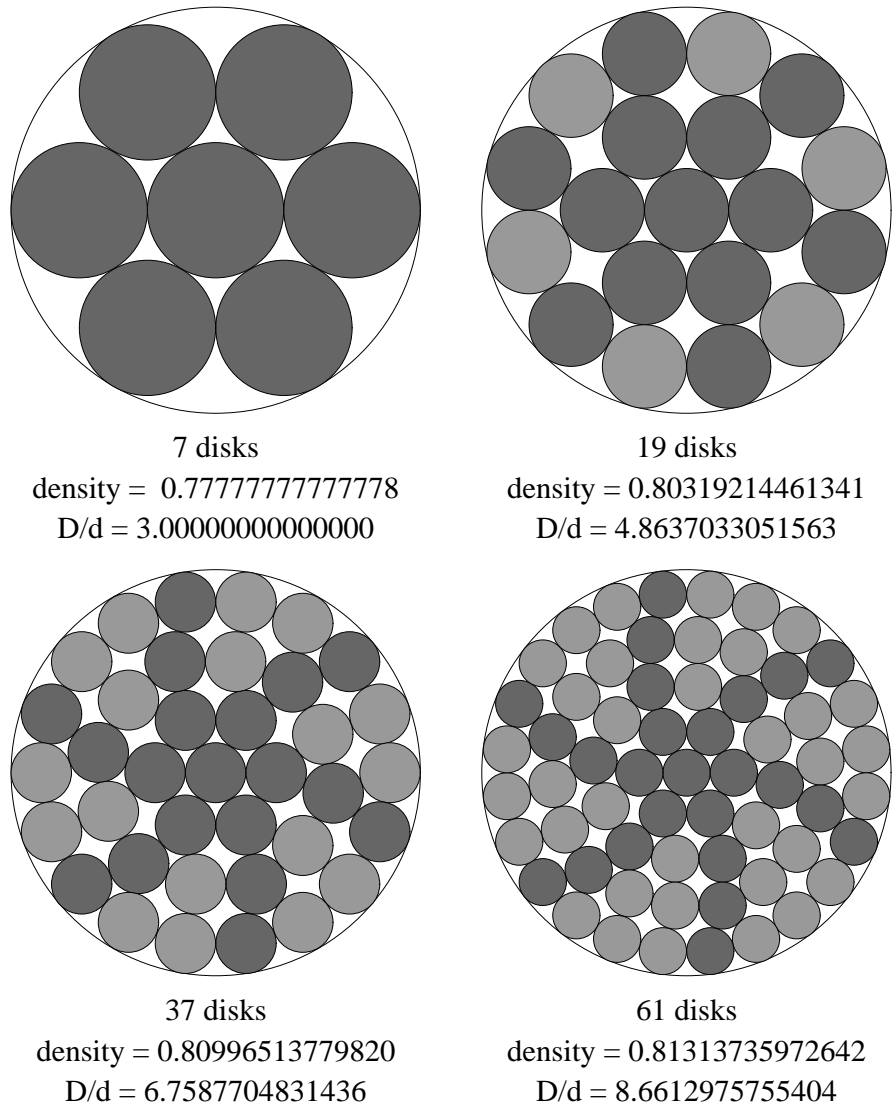,width=4.795in}}
\caption{ {\it Top}. The well known best packing of $h(1)=7$ disks
and the best previously conjectured packing of $h(2)=19$ disks.
{\it Bottom}. The best (that we found) packing of $h(3)=37$ disks and 
the best (one of the three that we found) packing of $h(4)=61$ disks.}
\label{4pack}
\end{figure}

\hspace*{\parindent}
The four best packings are presented in Fig.~\ref{4pack}.
The packing of 7 disks is well known and is easily seen 
to be optimal\footnote{
Indeed, if we partition the circle (of radius $r$) 
which contains the centers of the disks,
into six equal sectors, then
one sector must contain at least
two disk centers.
The distance between the two centers is at most
$r$, and the maximum is achieved when either one disk is centered
exactly at the enclosing circle center
or when two disks are centered at the peripheral corners of the sector.
If needed, we rotate the partition 
to exclude the latter possibility.
Then the central disk position defines the rest of
the packing pattern.}
and that of 19 disks is conjectured as such in \cite{K}.
The packings shown of 37 and 61 disks have not been reported before;
they are the best we found for these numbers of disks.
The density 7/9 of the 7-disk packing is presented
in decimal form in conformance with the other three densities;
an alternative finite form of the parameters
for all the packings in Fig.~\ref{4pack}
exists and is discussed in Section~3.
\section{The curved hexagonal pattern}
\hspace*{\parindent}
The pattern can be explained by comparing
it with the corresponding hexagonal pattern.
Fig.~\ref{twolarge} depicts the two patterns side-by-side
for $h(5)=91$ disks.
Each is composed of six sections which we tried to
emphasize by shading more heavily
the disks on the boundaries between them.
In the true hexagonal packing,
the sections are triangular,
the boundaries are six straight
``paths'' that lead from the central disk to
the six extreme disks.
In the curved hexagonal packing,
the triangular sections are ``curved'' 
and so are the paths.

To define the entire structure of the curved hexagonal packing,
it suffices to define positions of disks on one path.
Let the central disk be labeled 0 and the following disks
on the path be successively labeled $1,~2,~\ldots~k$.
Consider $k$ straight segments connecting the centers
of the adjacent labeled disks: 
0 to 1, 1 to 2,$~\ldots~,(k-1)$ to $k$.
Given a direction of rotation (it is clockwise in Fig.~\ref{twolarge}),
each following segment is rotated at the same angle, 
let us call it $\alpha$,
in this direction with respect to the previous segment.
Consider $\alpha$ as a parameter.
When $\alpha = 0$ we have the original hexagonal packing at the left.
If we gradually increase the $\alpha$,
``humps'' grow on the six sides of the hexagon.
When $\alpha$ reaches the value $\alpha_k = \frac {\pi} {3k}$,
all the disks of the last, $k$th, layer 
are at the same distance from the central disk.

\begin{figure}[htb]
\centerline{\psfig{file=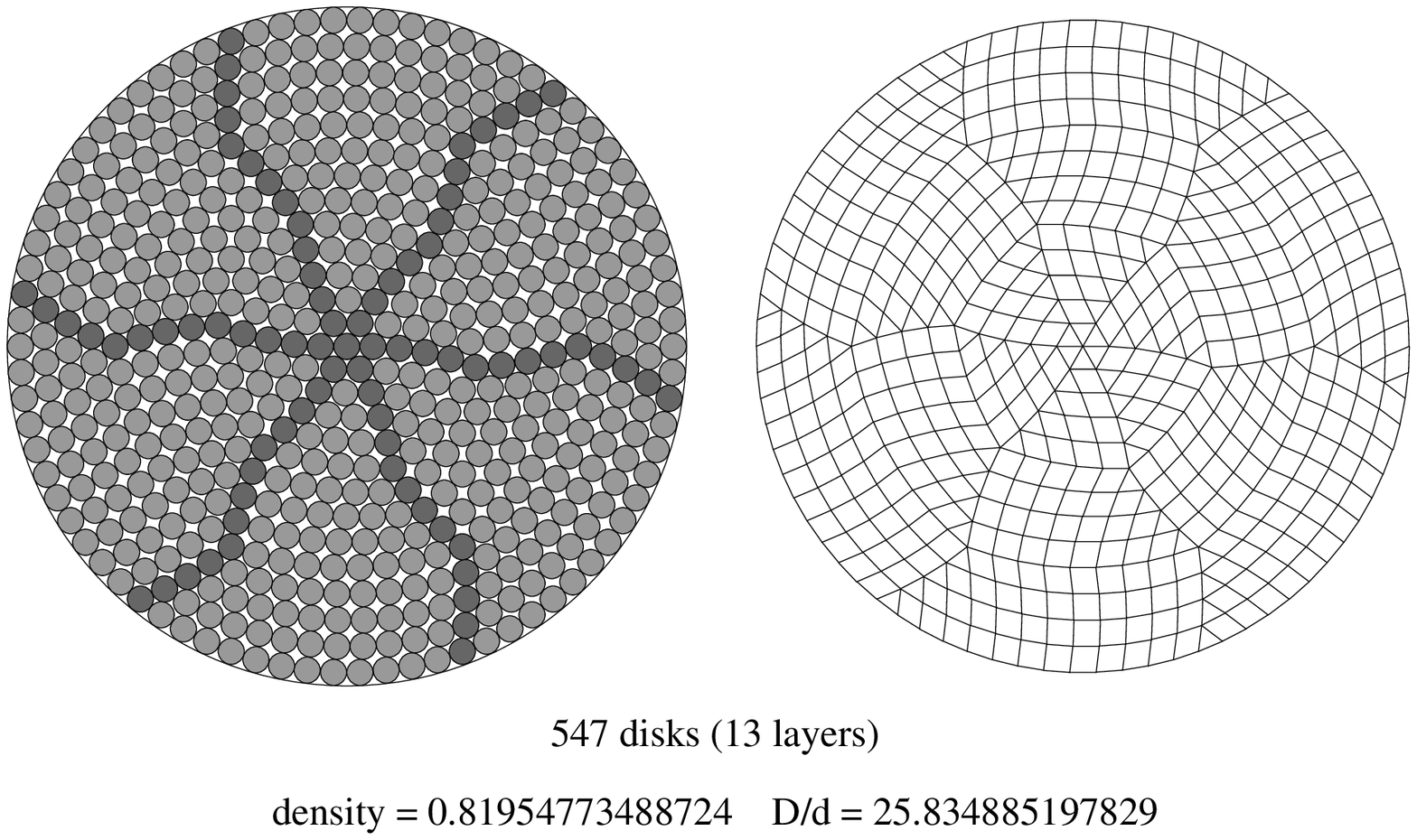,width=6.7in}}
\caption{{\it Left}. A curved
hexagonal packing of $h(13)=547$ disks in a circle.
In this packing the sense of rotation is flipped for
layers 6,7,8, and 9.
{\it Right}. The bond diagram for the packing at the left.
The diagram contains a straight segment
between centers of any pair of disks that are in contact.
}
\label{13layers}
\end{figure}

Following the path of labeled disks as defined above,
the distance $P$ in disk diameters ($d$)
from the center of disk 0 to the center of disk $k$ is given by
\beql{eq1}
P = | 1 + e^{i \alpha_k} + e^{2i \alpha_k} +  ~\ldots~ + e^{(k-1)i \alpha_k} | 
\eeq
Since $\alpha_k = \frac {\pi} {3k}$,
this simplifies to
\beql{eq2}
P =  \frac {1} {2 \sin \frac {\pi} {6k} }
\eeq
It follows that
the ratio $D/d$ of the enclosing circle diameter 
to the disk diameter for the curved hexagonal packing is 
\beql{eq3}
\frac {D} {d} = 1 + \frac {1} {\sin {\frac {\pi} {6k}}}
\eeq
The packing density,
i.e., the fraction of the enclosing circle area
which is covered by disks,
can then be found as
\beql{eq4}
density = h(k) (d/D)^2
\eeq
The density tends to the limit 
$\frac {\pi^2} {12} = 0.822467033...$ as $k \rightarrow \infty$.

The limit density of the curved hexagonal packing
is exactly the square of 
$\frac {\pi} {2 \sqrt{3}} = 0.906899682...$,
the density of the hexagonal packing of
the infinite plane (cf. \cite{FEJ}, \cite{FG}, \cite{GL1} or \cite{O}).
The fact that the latter density is larger than the former
implies that the curved hexagonal pattern is non-optimal for large $k$.
This is so, because as $n$ increases the density of 
hexagonal configurations of $n$ congruent disks
that fit inside a circle 
arbitrarily closely approximates $\frac {\pi} {2 \sqrt {3}}$
(see \cite{FEJ}).
Indeed, we found better packings for $k = 6$, 7, and 8
which we discuss in Section 5.

A curved hexagonal packing of $h(k)$ disks 
can be constructed for any $k \ge 1$. 
Fig.~\ref{13layers} depicts an instance for $k=13$.
We believe there are total of $m(k) = \max \{ \frac {(k-1)!} {2},1 \}$
non-congruent equally good
curved hexagonal patterns of $h(k)$ disks.
A curved hexagonal packing pattern is invariant under
a $60 ^{\circ}$ rotation.
A method to generate different curved hexagonal packings of 
a $k$-layered pattern, that is, for $h(k)$ disks, is as follows.
Take the pattern described above --
let us call it {\em basic} pattern for brevity --
and choose a subset among layers $2,~3,~\ldots,~k-1$.
Flip the sense of rotation of the chosen layers 
(or of the corresponding segments on the labeled path).
The flip is equivalent to making the mirror reflection of these layers.
The flipped layers will fit in their place because 
of mirror symmetry of segments of the layers enclosed between 
consecutive paths in the basic pattern.
Note that layers 1 and $k$ are not subject to the flip,
because the disks are positioned on them invariantly for all 
curved hexagonal packings for the given $k$.
The resulting $2^{k-2}$ combinations yield
$2^{k-3}$ different packings considering the reflection symmetry.
We call these modifications ``regular'' curved hexagonal packings
and Fig.~\ref{13layers} displays one of them.
\begin{figure}[htb]
\centerline{\psfig{file=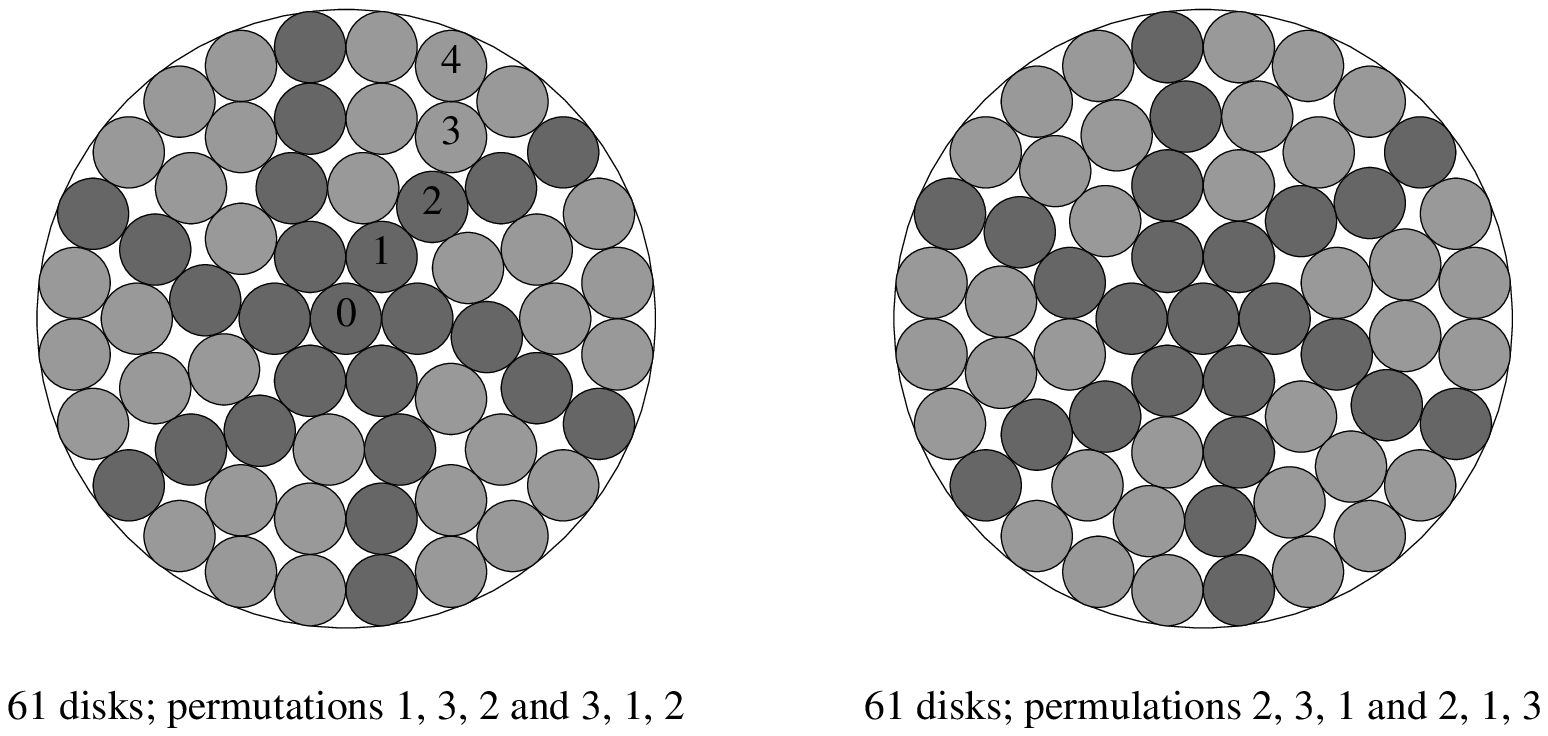,width=7in}}
\caption{Two out of the three existing 
packings of 61 disks in a circle; 
the third packing which corresponds to permutations 
1, 2, 3 and 3, 2, 1
is given in Fig.~2.1. The role of permutations is explained
in the text. The packing at the left is regular which is demonstrated
by the path of disks labeled 0, 1, 2, 3, 4,
each of which has a triangular hole attached. 
The packing at the right is irregular.}
\label{two61}
\end{figure}

We separate regular packings into a distinctive class 
because their existence 
follows from the existence of the basic regular packing
and it is easy to understand.
Besides, the regular packings are the most frequent ones among the optimal 
packings spontaneously generated by our ``billiards'' procedure 
(see Section 4).
However, the regular packings do not exhaust all the variants since
$2^{k-3} < m(k)$ for $k \ge 4$.
We believe that the packings in the full set can be identified by different
permutations in the order of summands 
$e^{i\alpha_k}, e^{2i\alpha_k},\ldots,e^{(k-1)i\alpha_k}$
in the expression (\ref{eq1}) for $P$,
or simply by permutations in the sequence
\beql{eq5}
1, 2, \ldots , k-1.
\eeq
For each permuted sequence, a path of $k+1$ disks can be constructed
and this path, when completed with layers, happens
to form a curved hexagonal packing which is
equal in quality to the basic one.
Permutation $i_1, \ldots , i_{k-1}$ of sequence (\ref{eq5})
produces the mirror reflection to the
pattern produced by permutation $k-i_1, \ldots , k-i_{k-1}$.

The method described above fills the pattern layer-by-layer
beginning from the central disk out.
An alternative method \footnote {Suggested by the anonymous referee}
fills the layers in the opposite direction.
The outermost layer contains $6k$ disks 
densely placed at the periphery.
(Examining the outermost layer is an alternative
way to infer (\ref{eq2}).)
Once layer $k$ is in place,
we choose a spot to attach the first disk of layer $k-1$
so that the disk contacts two disks of layer $k$.
There are $6k$ different attachment spots.
We attach the second, the third, ... the $6(k-1)$th
disk, say, clockwise, so that each next disk contacts the previous
disk and at least one disk of layer $k$.
As a result exactly 6 disks of layer $k-1$,
namely, disk 1, disk $k$, disk $2k-1$, disk $3k-2$, disk $4k-3$,
and disk $5k-4$ will each contact
two disks of layer $k$, while the remaining $6(k-2)$ disks
will each contact only one disk of layer $k$.
All thus obtained two-layer  configurations 
are congruent to each other by rotation.
However, beginning with layer $k-2$ (when $k \ge 3$), 
as we fill the pattern inward,
the choice of the spot for the placement of the ``first'' disk 
(the one that contacts two disks of the previous layer) distinguishes
the pattern, modulo a $60 ^{\circ}$ rotation.
Thus, we have $k-1$ different ways to place the first disk in layer
$k-2$, then $k-2$ ways to place the first disk in layer $k-3$, and so on.
This process yields $(k-1)!$ different clockwise placements.
Congruence by the reflection symmetry makes us to half the total: 
we have $\frac {(k-1)!} {2}$ non-congruent placements (for $k \ge 3$).

One can enumerate curved hexagonal packings
either by the sequence of rotation angles on a path
beginning from the central disk outward or by the sequence of
relative positions of ``first'' disks beginning
from the external layer inward.
For either method a simple computer program can be 
written that synthesizes a curved hexagonal 
packing given the sequence as an input.
In what follows we adopt the former enumeration
because the program we wrote uses that method.
Note that we do not offer here a formal proof that either
construction method described above actually works,
i.e., produces $\frac {(k-1)!} {2}$ different packings (for $k \ge 3$), 
as claimed, but we strongly believe it does.

The basic packing of 61 disks depicted in Fig.\ref{4pack} 
corresponds to the original sequence 1, 2, 3 or its reflection 3, 2, 1.
The packings of 61 disks that correspond to permuted sequences
1, 3, 2 and 2, 3, 1 or their reflections 3, 1, 2 and 2, 1, 3, respectively,
are shown in Fig.~\ref{two61}.
A characteristic feature of a regular packing is the existence
of a path any disk on which has a triangular hole attached.
On the {\em bond diagram} which has a straight segment
between centers of any pair of disks that are in contact
(see an example in Fig.~\ref{13layers}),
this path is seen
as a (broken) line of $k$ segments that leads
from the center to the
periphery each segment on which is a side of a triangle.
We identified
such a path in the packing at the left in Fig.~\ref{two61}
by labeling its disk components.
Hence this packing is regular.
The one at the right is irregular,
because no required path can be found for it.
The basic packing of 61 disks shown in Fig.~\ref{4pack}
is also regular.
This amounts to three different packings of $h(4)=61$ disks,
two packings of which are regular.
Indeed, $m(4) = 3$ and $2 = 2^{4-3}$.
Our belief is that for $n=61$ there is no better packing
than the three presented.
We also believe that there is no {\em fourth} packing 
equal in quality to the three presented
but distinct from any of them.

\begin{figure}[H]
\centerline{\psfig{file=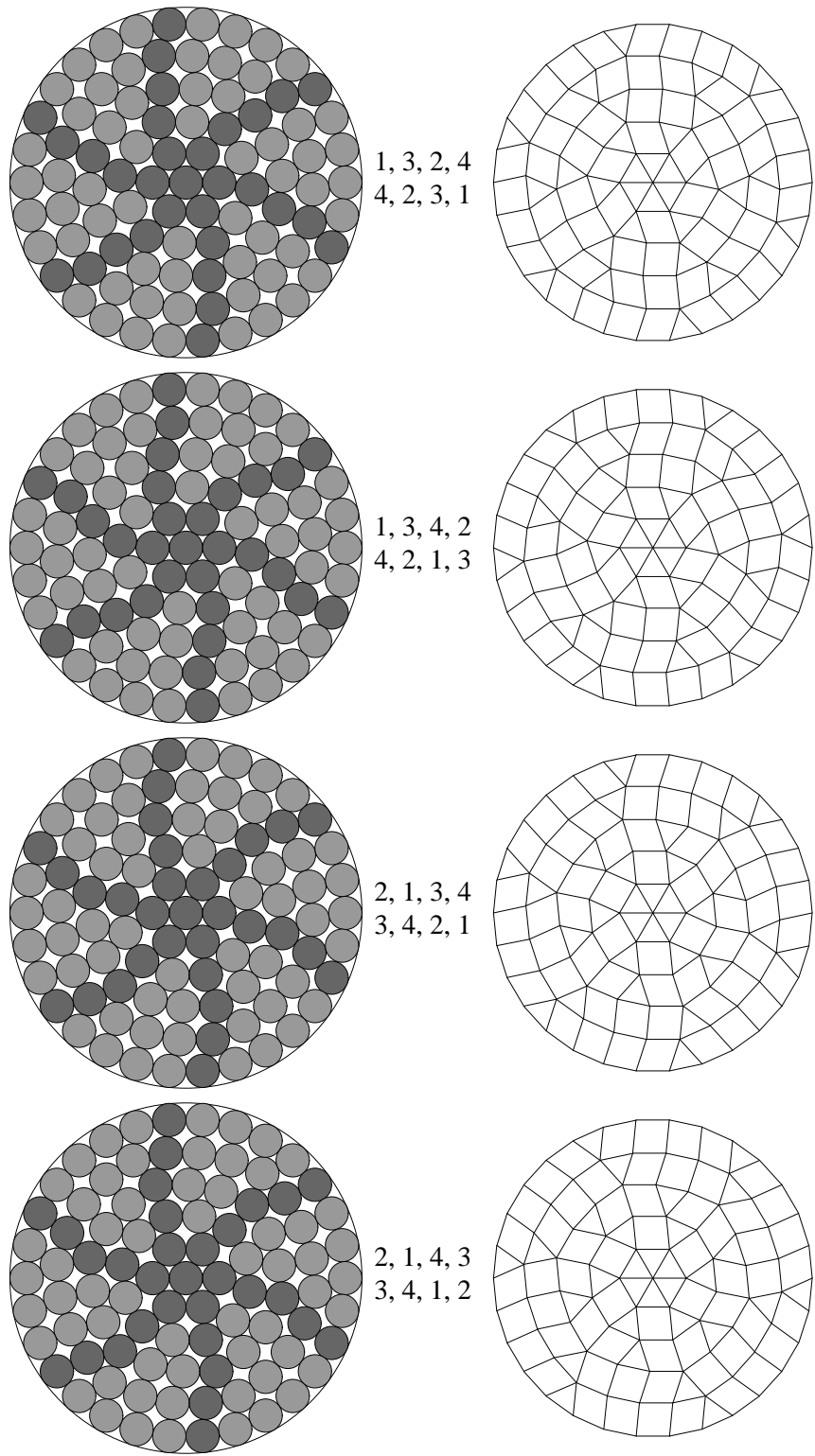,width=4.7in}}
\caption{First four irregular densest (that we found) 
packings of 91 disks.
Each packing shown at the left is accompanied by its two generating
sequences in the middle and the bond diagram 
at the right.}
\label{all91a}
\end{figure}

\begin{figure}[H]
\centerline{\psfig{file=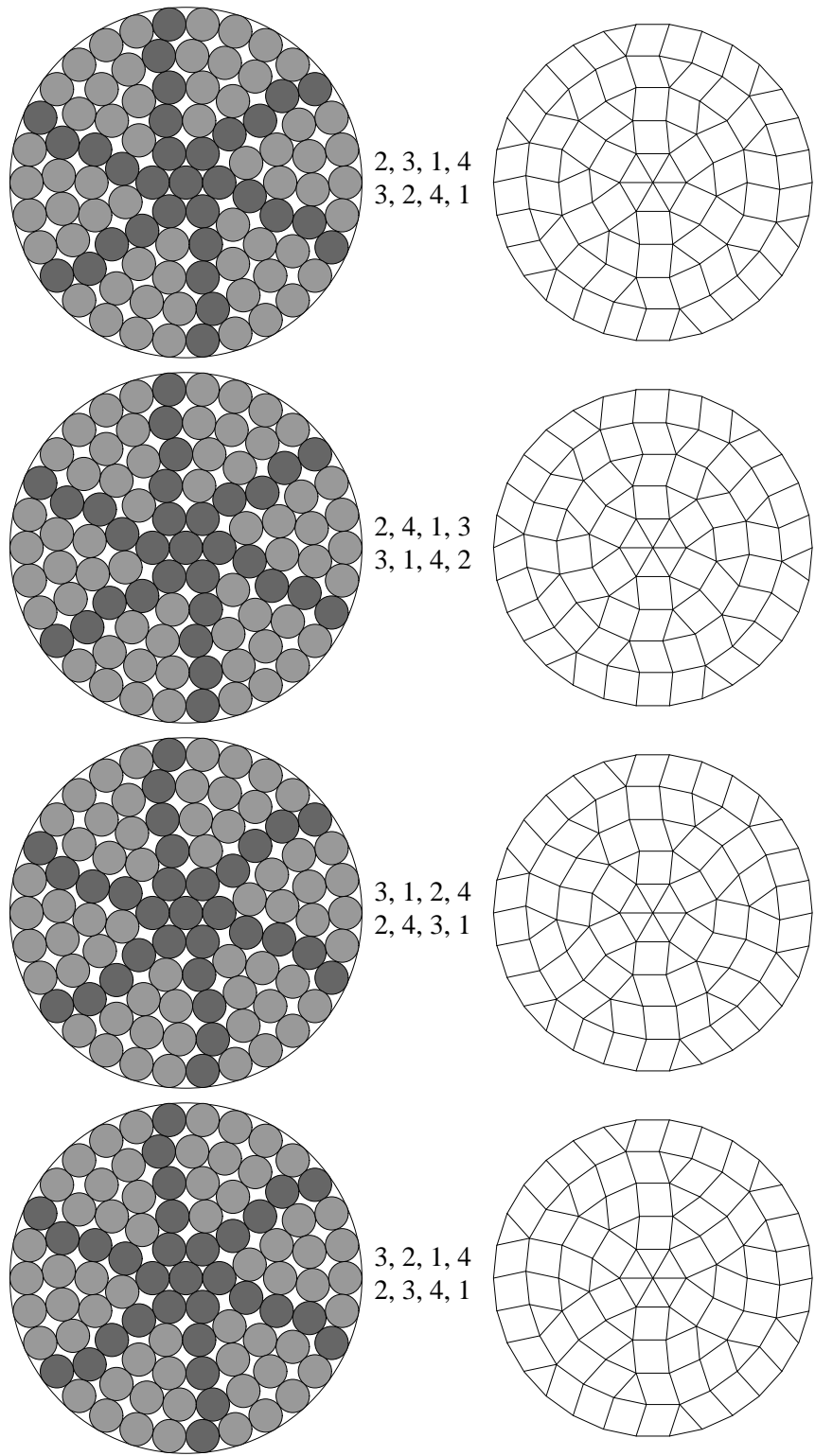,width=4.7in}}
\caption{Second four irregular densest (that we found)
packings of 91 disks.
Each packing shown at the left is accompanied by its two generating
sequences in the middle and the bond diagram 
at the right.}
\label{all91b}
\end{figure}

Similarly, there are
$m(5) = 12$ curved hexagonal packings of $h(5)=91$ disks,
$4 = 2^{5-3}$ packings of which are regular.
One regular packing is represented in Fig.~\ref{twolarge};
it corresponds to the original sequence 1, 2, 3, 4 
or its reflection 4, 3, 2, 1.
The three other regular packings 
correspond to permuted sequences (1, 2, 4, 3), (1, 4, 2, 3),
and (1, 4, 3, 2) or their reflection, respectively,
(4, 3, 1, 2), (4, 1, 3, 2), and (4, 1, 2, 3).
While the existence of the additional regular packings is obvious
and those packings are not shown,
for a reader who might be not fully convinced in the formal validity
of the construction procedures described above we show
eight irregular packings in Figs.~\ref{all91a} and \ref{all91b}.
Each packing is accompanied by its two generating sequences
and its bond diagram.
The bond diagrams clearly distinguish the packings 
and also show that the packings are irregular.
Again, we believe there is no packing better than those 12
and there is no 13th packing equal in quality to those 12
but distinct from any one of them.

For $k > 5$ there are packings of $h(k)$ disks
that are better than the $m(k)$ curved hexagonal ones.
It is therefore
easy to produce for such a $k$ infinitely many different configurations
of non-overlapping disks that are equal in quality to curved hexagonal ones.
We believe, however, that if a configuration of $h(k)$ congruent disks 
of the same quality as a curved hexagonal packing of the same disks 
has a general structure of a curved hexagonal packing described below,
then it must be one of the $m(k)$ curved hexagonal packings. 
The general structure consists of $k$
layers that surround the central disks;
layer $i$ consists of $6i$ disks placed in a circular fashion,
for $i=1,2,...k$, a disk of layer $i$ can only contact 
two disks in the same layer (previous and next disk along the circle) and
disks of layer $i-1$ and $i+1$ (with the central disk 
being counted as layer 0 and the boundary as layer $k+1$).

Note that we distinguish between a {\em configuration}
and a ({\em rigid}) {\em packing}.  We call a {\em configuration}
of non-overlapping disks inside a circle 
a {\em rigid} {\em packing} (or simply a {\em packing})
if there exists a non-empty subset of disks
in the configuration such that the only continuous motion of some or all
disks in the subset is rotation of this subset as a whole with
the center of rotation being at the center of the enclosing
container-circle.
Curved hexagonal patterns are rigid packings
because each internal $i$th layer of disks has triangles 
(as seen on bond diagrams)
that connect it to the corresponding outer layer $i+1$,
and because the outermost $k$th layer is obviously rigid.

\section{How the ``billiards'' algorithm produces packings}
\begin{figure}[H]
\centerline{\psfig{file=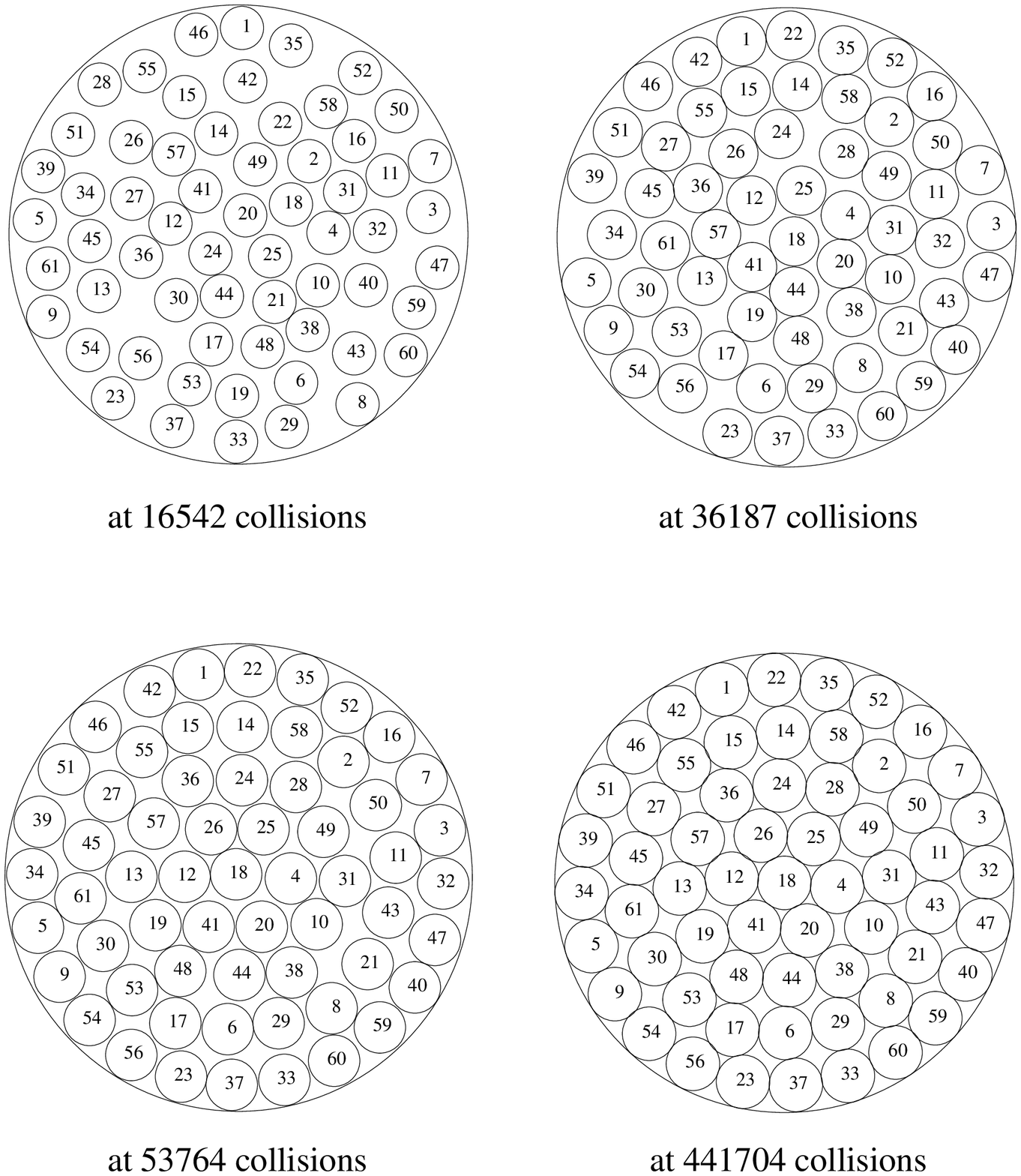,width=7in}}
\caption{Successive snapshots of simulating expansion
of 61 disks inside a circle.
The progress is monitored by counting collisions.
Disks are labeled 1 to 61 arbitrarily
but the same disk carries the same label in all four snapshots.
The last pattern is a mirror-reflection
of the packing of 61 disks in Fig.~2.1
}
\label{progress}
\end{figure}

\hspace*{\parindent}
A detailed description of the philosophy,
implementation and applications of this event-driven
algorithm can be found in \cite{L}, \cite{LS}.
Essentially, the algorithm simulates a system of $n$ perfectly elastic disks.
In the absence of gravitation and friction,
the disks move along straight lines,
colliding with each other and the region walls
according to the standard laws of mechanics,
all the time maintaining a condition of no overlap.
To form a packing,
the disks are uniformly allowed to gradually increase in size,
until no significant growth can occur.
Fig.~\ref{progress} displays four successive snapshots in an experiment
with 61 disks. 
We took the snapshots beginning the time 
when a local order begins to form and till the time
when the set of neighbors of each disk is stabilized.

The latest snapshot shown in Fig.~\ref{progress},
the one at 441704 collisions, looks dense
and, within the drawing resolution,
it is identical (up to a mirror-reflection) to
the final snapshot (see Fig.~\ref{4pack}).
However, numerically there are gaps
of the order of $10^{-5}$ to $10^{-3}$ of the disk diameter
in disk-disk and disk-wall pairs
that appear to be in contact with each other in Fig.~\ref{progress}.
Accordingly, only the first 3 decimal digits of $D/d$ corresponding 
to the latest shown in Fig.~\ref{progress}
snapshot are identical with the correct $D/d$.
To close these gaps and to achieve full convergence
of $D/d$ and of the density,
it usually takes 10 - 20 million further collisions.
We consider $D/d$ and the density to have converged when their values
do not change with full double precision for several million
collisions. 

Of course, 
as is typical in numerical iterative convergent procedures,
if the computations were performed with
the infinite precision, 
the convergence would be never achieved
and the ever diminishing gaps would always be there.
The ``experimental'' converged values agree
with the ``theoretical'' ones computed by formulas 
(\ref{eq3}) and (\ref{eq4}) to 14 or more significant digits.
Moreover, when we initialize the disk positions
differently, then the final parameters achieved
are either quite distinct and significantly smaller than those achieved
in the run presented in Fig.~\ref{progress} -- and then the corresponding
pattern is different from a curved hexagonal packing --
or they are identical to 14 or more significant digits
-- and then the corresponding final pattern is one of
the six known curved hexagonal ones.
This makes us suspect that we have found
the best possible packing and that its parameters 
$D/d$ and density are correct to 14 significant digits.

This algorithm does not equally favor 
the existing curved hexagonal packings.
For the chosen algorithm parameter settings,
including the slow disk expansion
(the ratio of disk expansion speed to the average
linear motion speed is 0.001),
the overwhelming majority
of produced curved hexagonal packings
were of the regular patterns.
For example,
out of our 1002 runs with $n = h(5) = 91$ disks,
curved hexagonal packings were obtained 90 times (9\%). 
Among those,
the four existing regular patterns (33\%) 
were seen in 81 runs (90\%),
with about the same frequency each.
Only three out of the existing eight irregular patterns
were seen in the remaining 9 runs.
We gave up waiting for the other five irregular patterns,
shown in Figs.~\ref{all91a} and \ref{all91b},
to be generated spontaneously, i.e., from random initial configurations.
Instead, we constructed those 
from their path sequences using the method discussed in Section~3.
\section{Packings of 127, 169, and 217 disks}
\begin{figure}[H]
\centerline{\psfig{file=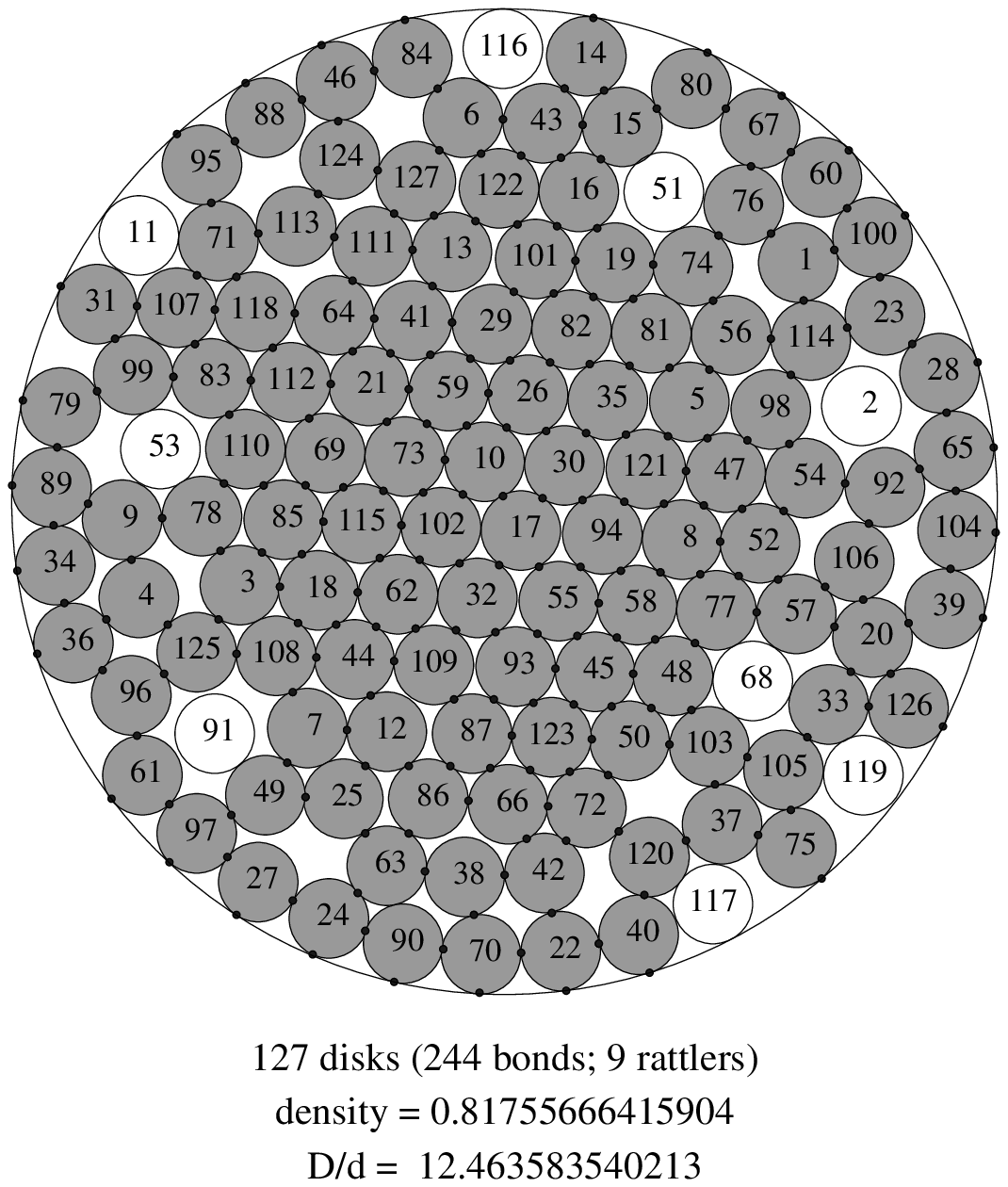,width=5in}}
\caption{A packing of $h(6)=127$ disks in a circle
that is better than the corresponding curved hexagonal packings.
Each disk is provided with its unique identification label
to facilitate the reference.
Little black dots are ``bonds''; a bond indicates that
the corresponding distance is less 
than $10^{-13}$ of the disk diameter.
Where a pair disk-disk or disk-wall are apparently in contact
but no bond is shown, e.g., between disks 3 and 108,
the computed distance is at least $10^{-5}$ 
of the disk diameter.
The shaded disks can not move given the positions of their neighbors,
the non-shaded are ``rattlers'' that are free to move within their confines.
It might seem that there are sufficiently large cavities near 
disks 27, 104, 28, 88, that are in contact with the boundary, so 
that if the disks are pushed into the cavities 
the packing will ``unjam.'' This does not happen.
For example, pushing disk 27 into the position 
of contact with 97 and 49 results in the overlap with 24.
}
\label{127disk}
\end{figure}

\hspace*{\parindent}
We ran the ``billiards'' algorithm for
$n = h(6) = 127$, $n = h(7) = 169$, and $n = h(8) = 217$ disks.
For these $n$ the algorithm produced better packings 
than the curved hexagonal ones.
The patterns of those packings can all be described
as a (possibly disturbed) hexagonal disk assembly in the middle
surrounded with
irregularly placed disks at the circular border.
For larger $n$ this common pattern becomes more evident.
As an example, we show in Fig.~\ref{127disk}
the best packing we obtained for 127 disks.

We believe the packings achieved for $n = 127$, 169, and 217 are stable.
But we do not think they are the best,
because of a large number of local minima for these $n$.
For example, the packing shown in Fig.~\ref{127disk}
was the best among 111 independent tries.
Each try resulted in a packing
which is distinct from the others and had distinct
parameters $D/d$ and density.
14 out of 111 packings were better than the curved hexagonal.

These runs for $h(k)$ disks, $k > 5$,
were in contrast to the runs for $n = h(5) = 91$
and smaller $n$.
For example, for 91 disks the best (curved hexagonal) packings
were obtained 90 times out of 1002
with parameters $D/d$ and density
agreeing to 15 significant places.
The results of packings for $h(k)$, $k > 5$, disks are summarized in
Table~\ref{summary}.
\begin{table}[H]
\begin{center}
\caption{Curved hexagonal packings vs. the experimental ones 
for $k=$6, 7, and 8}
\label{summary}
~ \\
\begin{tabular}{| c |  r | r | r |}  \hline
 $k$    &    6   &    7    &     8    \\     \hline  
 $n$    &   127 &   169   &     217  \\     \hline \hline
curv. hex. density      &0.81622935362082 &0.81710701192903  &0.81776562948873 \\   \hline
curv. hex. $D/d$        &12.473713245670&14.381489999655&16.289788298679 \\   \hline\hline
experimental density    &0.81755666415904&0.82672262170717&0.83499393075147 \\ \hline
experimental $D/d$      &12.463583540213&14.297609837687&16.120860041887 \\ \hline
better packings        &14 out of 111& all 70  &  all 62    \\ \hline
\end{tabular}
\end{center}
\end{table}

\section{Tightness of curved hexagonal packings}
\begin{figure}[H]
\centerline{\psfig{file=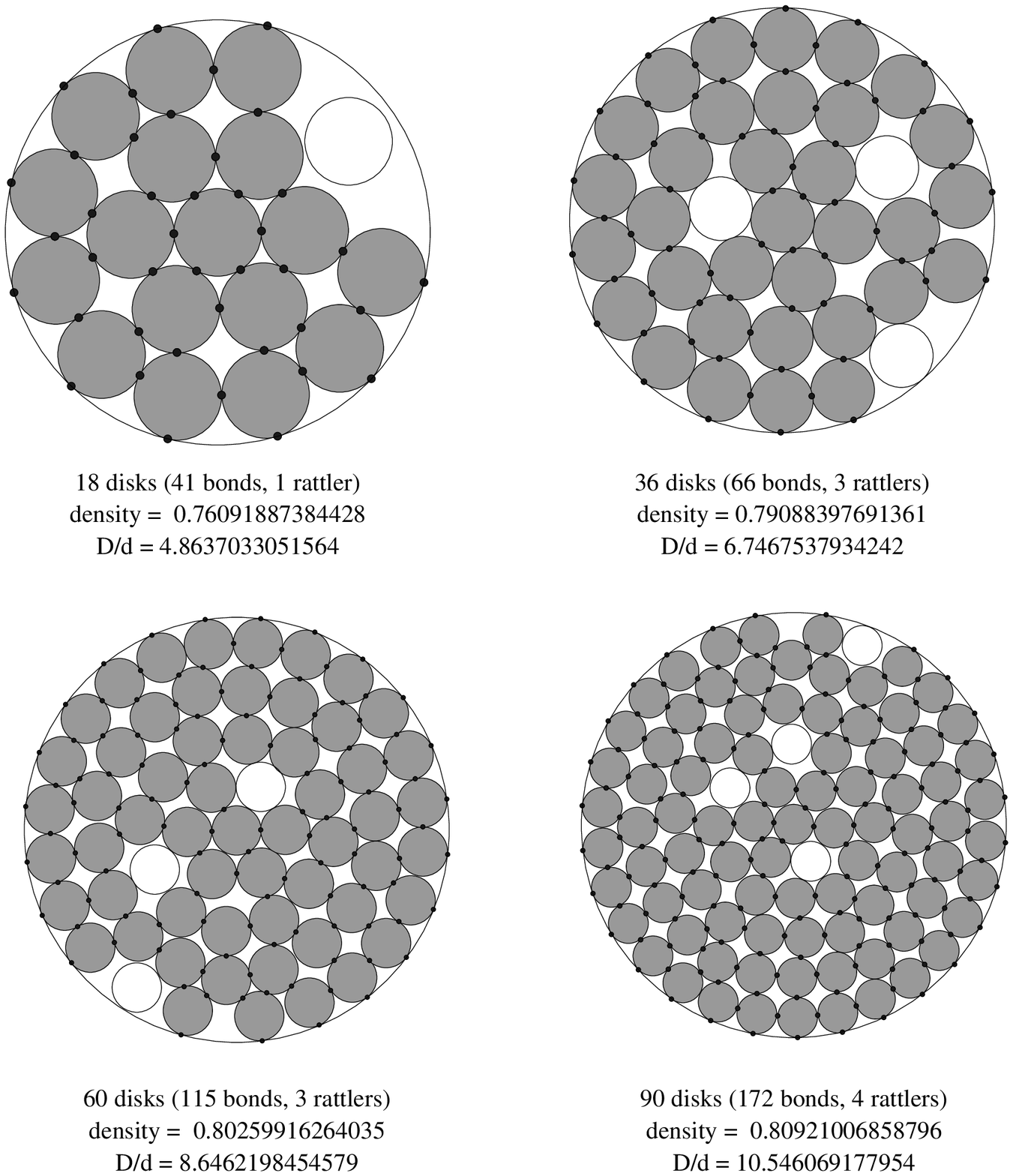,width=7in}}
\caption{The best found packings of $n=h(k)-1$ disks for
$k=2$, 3, 4, and 5 ($n=18$, 38, 60, and 90).}
\label{8pack4a}
\end{figure}

\begin{figure}[H]
\centerline{\psfig{file=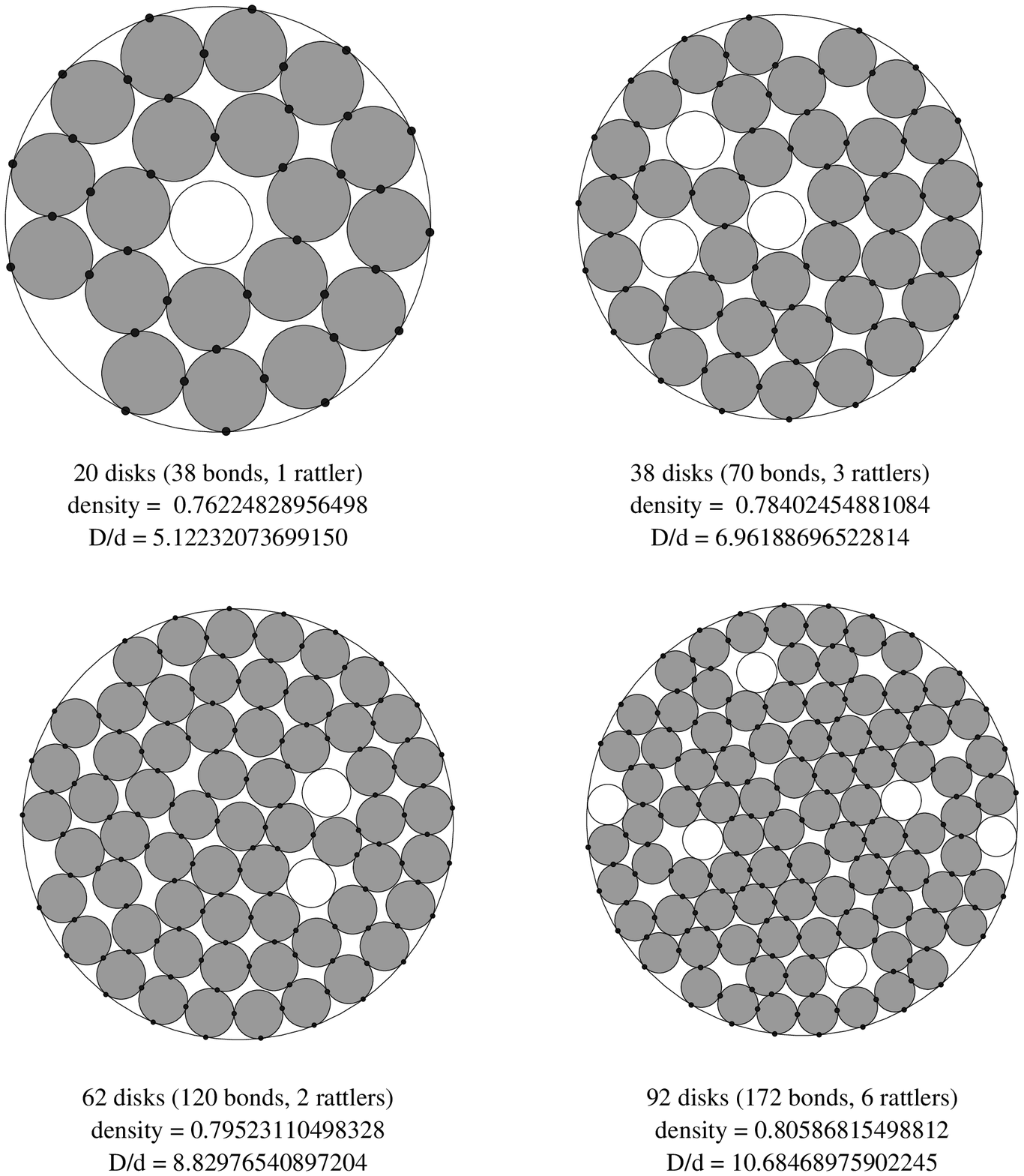,width=7in}}
\caption{The best found packings of $n=h(k)+1$ disks for
$k=2$, 3, 4, and 5 ($n=20$, 38, 62, and 92).}
\label{8pack4b}
\end{figure}

\hspace*{\parindent}
Fig.~\ref{8pack4a} and Fig.~\ref{8pack4b}
depict the best found packings of $n=h(k)-1$ and $h(k)+1$
disks for $k=2$, 3, 4, and 5,
that is, for $n=18$, 20, 36, 38, 60, 62, 90, and 92.
\footnote {The forthcoming paper \cite{GLNO} lists 9 more distinct
packings of 18 disks that are equal in quality to the one 
presented in Fig,~\ref{8pack4a}.}
The known packings of $n=6$ and 8 disks have
to the same tendencies, namely:

(a) the pattern of dense packing of $n = h(k)-1$ disks is obtained
by removing one disk from the pattern of dense packing of $h(k)$
disks (which is a curved hexagonal packing for the considered $k \le 5$)
and, for $k \ge 3$, by a small rearrangement of the disks;
its parameter $D/d$ is either not changed (for $k=1$ and $2$)
or is decreased only slightly (for $k=3$, 4, and 5);

(b) the pattern of dense packing of $n = h(k)+1$ disks differs
significantly from the pattern of dense packing of $h(k)$
disks and its parameter $D/d$ is increased
substantially for all $k=1$, 2, 3, 4, and 5).

The changes in $D/d$ are ``slight'' and ``substantial'' only
in comparison to each other.
The ratio of the decrease of $D/d$ for $n=h(k)-1$
over the increase of $D/d$ for $n=h(k)+1$
is given in Table~\ref{tight}.
\begin{table}[H]
\begin{center}
\caption{The ratio of the decrease of $D/d$ for  $n=h(k)-1$ 
over its increase for $n=h(k)+1$.}
\label{tight}
~ \\
\begin{tabular}{| c |  r | r | r | r | r |}  \hline
$k$       & 1    &   2    &    3    &   4    &   5  \\    \hline
$n=h(k)$  & 7    &  19    &   37    &  61    &  91  \\    \hline
ratio     & 0    &   0    &   0.0592 &  0.0895 &  0.1755  \\   \hline
\end{tabular}
\end{center}
\end{table}
\section{Discussion}
\hspace*{\parindent}
Our experiments reveal that for a sufficiently large $n$
good packings of $n$ equal disks in a circle have 
a complex pattern like that in Fig.~\ref{127disk},
with a large, perhaps disturbed, core of hexagonally packed disks
and irregularly placed disks along the periphery.
The fraction of the peripheral irregularity disks and the
perturbation in the hexagonally packed core usually diminish with $n$.
It is very difficult to obtain the best packings for large $n$
but we would guess their pattern to be of the same irregular type.
On the other hand, for $n \le 25$ 
symmetric and regular patterns of the best packing 
have been previously observed
that do not obey
the general description
given above for large $n$.
Our computer experiments show that at least for a particular
class of $n=h(k)=3k(k+1)+1$, the transition
from the regular, here curved hexagonal, pattern
to the irregular core-hexagonal one
occurs between $n=91$ and $n=127$.

\end{document}